\documentclass[11pt]{amsart}
\usepackage{amstext}
\usepackage{amsfonts}
\usepackage{amssymb}
\usepackage{amsbsy}
\usepackage{amsmath}
\usepackage{latexsym}
\usepackage{hhline}

\newcounter{num}[section] %

\newenvironment{theo}
{\refstepcounter{num}%
\bigskip\noindent{\bf Theorem~\arabic{section}.\arabic{num}. }\it}

\newenvironment{cor}
{\refstepcounter{num}%
\bigskip\noindent{\bf Corollary~\arabic{section}.\arabic{num}. }\it}

\newenvironment{lemma}
{\refstepcounter{num}%
\bigskip\noindent{\bf Lemma~\arabic{section}.\arabic{num}. }\it}

\newenvironment{example}
{\refstepcounter{num}%
\bigskip\noindent{\bf Example~\arabic{section}.\arabic{num}.}}

\newenvironment{remark}
{\refstepcounter{num}%
\bigskip\noindent{\bf Remark~\arabic{section}.\arabic{num}.}}

\newcounter{thepic}

\newcommand{\si}{\sigma}
\newcommand{\al}{\alpha}
\newcommand{\be}{\beta}
\newcommand{\ga}{\gamma}

\newcommand{\de}{\delta}

\newcommand{\LA}{\langle}
\newcommand{\RA}{\rangle}

\newcommand{\Char}{\mathop{\rm char}}

\newcommand{\FF}{{\mathbb{F}}}   

\newcommand{\ZZ}{{\mathbb{Z}}}   

\begin{document}
\title[Identities of sum of two PI-algebras]{Identities of sum of two PI-algebras\\ in the case of positive characteristic}%
\author{Ivan Kaygorodov}
\address{Ivan Kaygorodov\\
Universidade Federal do ABC, CMCC, Santo Andr\'{e}, SP, Brazil}
\email{kaygorodov.ivan@gmail.com}

\author{Artem Lopatin}
\address{Artem Lopatin\\
Sobolev Institute of Mathematics, Omsk Branch, SB RAS, Omsk, Russia;
Universidade Estadual de Campinas, IMECC, Campinas, Brazil}
\email{artem\_lopatin@yahoo.com}

\author{Yury Popov}
\address{Yury Popov\\
Sobolev Institute of Mathematics, SB RAS, Novosibirsk, Russia}
\email{yuri.ppv@gmail.com}

\begin{abstract}
 We consider the following question posted by K.I.~Beidar and A.V.~Mikhalev in 1995 for an associative ring $R=R_1+R_2$: is it true that if the subrings $R_1$ and $R_2$ satisfy polynomial identities, then $R$ also satisfies a polynomial identity? Over a field of positive characteristic we establish new conditions on $R_1$ and $R_2$ that guarantee a positive answer to the question. We find upper and low bounds on the degrees of identities of $R$.


\end{abstract}

\maketitle

\section{Introduction}\label{section_intro}

\subsection{Sum of rings}
We assume that $\FF$ is a field of arbitrary characteristic $p=\Char{\FF}\geq0$. All vector spaces, algebras and modules are over $\FF$ and all algebras are associative, except those mentioned in Subsection~\ref{subsec:1.4}. An algebra and a ring may not have a unity.

For a commutative ring $K$, denote by $K\LA X\RA$ the free $K$-module, freely generated by all non-empty products of letters $x_1,x_2,\dotsc$. Given a ring $R$ of characteristic $m\geq0$, a polynomial identity of $R$ is a non-zero element $f(x_1,\ldots,x_n)$ of $\ZZ_m\LA X\RA$ such that $f(r_1,\ldots,r_n)=0$ in $R$ for all $r_1,\ldots,r_n\in R$, where $\ZZ_m$ stands for $\ZZ/m\ZZ$ in case $m>0$ and $\ZZ_0=\ZZ$. Similarly we can define a polynomial identity of an $\FF$-algebra $A$ as an element of $\FF\LA X\RA$. A ring (an algebra, respectively) that satisfies a polynomial identity is called a PI-ring (PI-algebra, respectively). For short, polynomial identities are called identities.

It was shown by Regev~\cite{Regev71, Regev72} in~1971 that the tensor product of two PI-algebras is a PI-algebra.
A similar question can be considered for sums of rings. Namely, let a ring $R=R_1+R_2$ be the sum of two subrings $R_1$ and $R_2$ (i.e., every element of $R$ is equal to the sum $r_1+r_2$ for some $r_1\in R_1$ and $r_2\in R_2$). In 1995 K.I.~Beidar and A.V.~Mikhalev \cite{BeidarMikhalev95} posed the following question which is still open:%

\medskip%

\noindent{}{{\bf Beidar--Mikhalev's Problem.}\it\;  Is it true that if $R_1$ and $R_2$ satisfy polynomial identities, then $R=R_1+R_2$ also satisfies a polynomial identity?}
\medskip

\noindent{}The same question can also be asked for algebras over a field $\FF$.

\subsection{Known results}
A positive answer to Beidar--Mikhalev's Problem is known in many cases. O.H.~Kegel~\cite{Kegel62} has established that if $R_1$ and $R_2$ are nilpotent, then $R$ is also nilpotent. By a result of Yu.~Bahturin and A.~Giambruno~\cite{BahturinGiambruno94}, if $R_1$ and $R_2$ are  commutative rings, then $R$ satisfies the identity $[x,y][a,b]=0$, where $[x,y]=xy-yx$. In the aforementioned paper by K.I.~Beidar and A.V.~Mikhalev~\cite{BeidarMikhalev95} it has been shown that if $R_1$ and $R_2$ satisfy the identity $[x_1,x_2]\cdots [x_{2n-1},x_{2n}]=0$, then $R$ is a PI-ring. In~\cite{KepczykPuczylowski96} and~\cite{KepczykPuczylowski98} M.~K\c{e}pczyk and E.R.~Puczy\l{}owski have established that if $R_1$ satisfies the identity $x^n=0$ and $R_2$ is a PI-ring, then $R$ is also a PI-ring. Note that the generalization of Beidar--Mikhalev's Problem to the case of three summands has a negative solution since it was shown by L.A.~Bokut'~\cite{Bokut76} that every algebra can be embedded into a simple algebra, which is a sum of three nilpotent subalgebras.

Beidar--Mikhalev's Problem is known to have a positive solution when some additional conditions are imposed on products of elements of $R_1$ and $R_2$. For instance, it was shown by L.H.~Rowen~\cite{Rowen76} in 1976 that the following condition suffices: both $R_1$ and $R_2$ are left (or right) ideals of $R$. Then M.~K\c{e}pczyk and E.R.~Puczy\l{}owski~\cite{KepczykPuczylowski01} have established that it suffices to assume that $R_1$ is a left or right ideal. Moreover, B.~Felzenszwalb, A.~Giambruno and G.~Leal~\cite{FelzenszwalbGiambrunoLeal03} have extended this result even further by proving that the problem has a positive solution in case $(R_1R_2)^k\subset R_1$ or $(R_1R_2)^k\subset R_2$ for some $k>0$. They have also found the following upper bound on the degree $D'$ of a polynomial identity that holds in $R$:
$$D'\leq a^a+1\;\text{ for }\; a=8e(kd(d-1)-1)^2(d-1)^2,$$
where $R_1$ and $R_2$ satisfy identities of degree $d$, and $e$ is the base of the natural logarithm.

An $\FF$-algebra $B$ is said to \emph{almost} satisfy some property if there exists a two-sided ideal $I$ of $B$ of finite codimension that satisfies this property. In 2008 M.~K\c{e}pczyk~\cite{Kepczyk08} established that if $R_1$ and $R_2$ are almost nilpotent $\FF$-algebras, then $R$ is also an almost nilpotent $\FF$-algebra. In a recent paper~\cite{Kepczyk15} M.~K\c{e}pczyk has extended the previous result as follows. Consider two identities $f_1$ and $f_2$ such that Beidar--Mikhalev's Problem has a positive solution for all algebras $R_1$, $R_2$ satisfying $f_1$ and $f_2$, respectively. Then any algebra $A=A_1+A_2$ with subalgebras $A_1$ and $A_2$ almost satisfying identities $f_1$ and $f_2$, respectively, is a PI-algebra.

\subsection{New results}
In this paper we consider the case of positive characteristic of the field $\FF$. Our main result (Corollary~\ref{cor2}) is the following generalization of~\cite{FelzenszwalbGiambrunoLeal03}: Beidar--Mikhalev's Problem has a positive solution for an algebra $R=R_1+R_2$ if the subalgebras $R_1$ and $R_2$ almost satisfy the property that $(R_1R_2)^k$ is a subset of $R_1$ or $R_2$ for some $k>0$. Note that in Corollary~\ref{cor2} we actually do not have to require that a subalgebra of finite codimension is an ideal (see Theorem~\ref{theo_main} below). In Corollary~\ref{cor1} we improve the upper bound from~\cite{FelzenszwalbGiambrunoLeal03} on the degree of an identity of $R$ in terms of degrees of the symmetric identities that hold in $R_1$ and $R_2$. Corollary~\ref{cor2} and Corollary~\ref{cor1} are partial cases of Theorem~\ref{theo_main}, which states that under certain conditions an identity of degree less than $C(d_1+d_2)d_1$ or $C(d_1+d_2)d_2$ holds in $R$, where $d_i$ is the minimal degree of the symmetric identity that holds in $R_i$ for $i=1,2$, and $C$ does not depend on $d_1$ and $d_2$. On the other hand, in general $R$ does not satisfy an identity of degree less than $d_1 d_2$, where $d_i$ is the least low degree of an identity that holds in $R_i$ for $i=1,2$, see Theorem~\ref{theo_low_bound}. (The definition of the least low degree is given in Section~3.) The last result holds over a field of arbitrary characteristic.

\subsection{Non-associative case}
\label{subsec:1.4}
An analogue of Beidar--Mikhalev's Problem can be considered for non-associative  algebras. In~\cite{Kegel62} O.~Kegel has asked whether a Lie ring that can be written as a sum of two nilpotent subrings is solvable. In the case of finite dimensional Lie algebras over a field this problem has been solved. Namely, the answer is positive for a field of characteristic 0  (see~\cite{Goto62}), as well as for a field of odd characteristic $p$. The last case has been independently considered by V.~Panyukov~\cite{Panyukov90} for $p>2$ and by  P.~Zusmanovich~\cite{Zusmanovich91} for $p>5$. On the other hand, A.~Petravchuk~\cite{Petravchuk88} has found a counter-example in the case of characteristic 2. Consider a solvable Lie algebra $L=L_1+L_2$ of the solvability degree $s$, where the nilpotency degrees of the subalgebras $L_1$ and $L_2$ are $n_1$ and $n_2$. The problem of describing an explicit upper bound on $s$ in terms of $n_1$ and $n_2$ is open. This problem is solved only in case $n_1=n_2=1$ (i.e., if $L_1$ and $L_2$ are abelian), where the answer is $s=2$ (see~\cite{Kolman65}).

Over a field of characteristic different from $2$ and $3$ S.~Pchelintsev established that any alternative algebra which is the sum of two solvable subalgebras is solvable itself.

\section{Symmetric identities}

In 1993 A.~Kemer~\cite{Kemer93} established that any PI-algebra over a field of positive characteristic satisfies the symmetric identity
$$s_d=\sum_{\si\in S_d} x_{\si(1)}\cdots x_{\si(d)}$$
for some $d$, where $S_d$ stands for the group of all permutations of $n$ elements. We say that a subset $L$ of an algebra $R$ satisfies an identity $f=f(x_1,\ldots,x_n)\in \FF\LA X\RA$ if $f(a_1,\ldots,a_n)=0$ for all $a_1,\ldots,a_n$ from $L$. The aforementioned result by Kemer implies that over a field of positive characteristic if a set satisfies an identity, then it satisfies some symmetric identity.

In the rest of this section we assume that $R=R_1+R_2$, where $R_1$ and $R_2$ are PI-algebras over a field $\FF$ of positive characteristic $p$, unless otherwise is stated.

\begin{theo}\label{theo_main} Let $R_1$ and $R_2$ satisfy the symmetric identities of degrees $d_1$ and $d_2$, respectively. Assume that there exist $k>0$ and subalgebras $L_i\subset R_i$ of finite codimensions $t_i$  ($i=1,2$) such that $(L_1L_2)^k$ is a subset of $R_i$ for $i=1$ or $i=2$. Then $R$ satisfies the symmetric identity of degree
$$D=(d_1+d_2-2)(k d_i + 2)((t_1+t_2)(p-1)+1).$$
\end{theo}

We present the proof of this theorem below.

\begin{cor}\label{cor2} Assume that $R_1$ and $R_2$ almost satisfy the following property: $(R_1R_2)^k$ is a subset of $R_1$ or $R_2$ for some $k>0$. Then $R$ is a PI-algebra.
\end{cor}

\begin{cor}\label{cor1} Let $R_1$ and $R_2$ satisfy the symmetric identities of degrees $d_1$ and $d_2$, respectively, and $(R_1R_2)^k$ is a subset of $R_i$ for $i=1$ or $i=2$. Then $R$ satisfies the symmetric identity of degree  $D=(d_1+d_2-2)(k d_i + 2)$.
\end{cor}

\medskip

To prove Theorem~\ref{theo_main} we will apply the following two lemmas.

\begin{lemma}\label{lemma1}
Suppose that $a_1,\ldots,a_d$, where $d>0$, are elements of an algebra $A$. Then $s_d(a_1,\ldots,a_d)=0$ whenever there exists $l\leq d$ such that $a_1,\ldots,a_l$ lie in some subspace $U$ of $A$ of dimension $t>0$ and $l> t(p-1)$.
\end{lemma}
\begin{proof}
Since $s_d$ is linear, it is enough to prove the lemma in case $a_1,\ldots,a_l$ are elements of some basis $\{u_1,\ldots,u_t\}$ of $U$. Moreover, for some $i_0$ a basis element $u_{i_0}$ appears in $\{a_1,\ldots,a_l\}$ at least $p$ times. Without loss of generality we may assume that $a_1=\cdots=a_p=u_{i_0}$.
Since
$$s_d(\underbrace{a_1,\ldots,a_1}_p,a_{p+1},\ldots,a_d)= p! \sum a_{\si(1)}\cdots a_{\si(d)}=0,$$
where the sum ranges over all permutations $\si\in S_d$ with $\si(1)<\cdots<\si(p)$, the lemma has been proved.
\end{proof}

\begin{lemma}\label{lemma_main} Assume that $\FF$ is a field of arbitrary characteristic. Let vector subspaces $V_1$ and $V_2$ of an algebra $R$  satisfy the symmetric identities of degrees $d_1$ and $d_2$, respectively. Assume that there exists $k>0$ such that the set $(\hat{V}_1\hat{V}_2)^k$ satisfies the symmetric identity of degree $d_3$, where $$\hat{V}_i=\{v_1\cdots v_r\,|\,  r<d_i \text{ and } v_1,\ldots,v_r\in V_i\}$$
for $i=1,2$. Then the subspace $V_1+V_2$ of $R$ satisfies the symmetric identity of degree
$$D=(d_1+d_2-2)(k d_3 + 2).$$
\end{lemma}
\begin{proof}
We denote the unity of $R$  by $1$ if this unity exists. Otherwise, we write $1$ for a formal element with the property $1\cdot a=a\cdot 1 = a$ for each $a$ from $R$.

We will show that $s_D(a_1,\ldots,a_D)=0$ in $R$ for all $a_1,\ldots,a_D$ from $V=V_1+V_2$.
Since $s_D$ is linear, without loss of generality we may assume that $a_1,\ldots,a_D$ are elements of $V_1\cup V_2$. Moreover,  we may assume that $a_1=b_1,\ldots,a_r=b_r$ are elements of $V_1$ and $a_{r+1}=c_1,\ldots,a_{D}=c_s$ are elements of $V_2$ for some $r,s\geq0$ satisfying the equality $r+s=D$. A sequence $\be=(\be_0,\be_1,\ldots,\be_m)$ with $m\geq1$ is called an $r$-partition if
$$0=\be_0\leq\be_1<\be_2<\cdots<\be_{m-2}<\be_{m-1}\leq \be_m=r.$$
Given $\si\in S_r$, an $r$-partition $\be=(\be_0,\ldots,\be_m)$ and $1\leq i\leq m$, let
$$b^{\be}_{\si}(i)=\left\{
\begin{array}{rl}
b_{\si(\be_{i-1}+1)}\cdots b_{\si(\be_i)}&\text{ if } \be_{i-1}<\be_i, \\
1,&\text{ if } \be_{i-1}=\be_i. \\
\end{array}
\right.$$
Similarly, for $\tau\in S_s$, an $s$-partition $\ga=(\ga_0,\ldots,\ga_m)$ and $1\leq i\leq m$, let
$$c^{\ga}_{\tau}(i)=\left\{
\begin{array}{rl}
c_{\tau(\ga_{i-1}+1)}\cdots c_{\tau(\ga_i)},&\text{ if } \ga_{i-1}<\ga_i, \\
1,&\text{ if } \ga_{i-1}=\ga_i. \\
\end{array}
\right.$$

Given $r,s\geq0$, denote by $\Delta_{r,s}$ the set of all  pairs $(\be,\ga)$ of an $r$-partition $\be=(\be_1,\ldots,\be_m)$ and an $s$-partition $\ga=(\ga_1,\ldots,\ga_m)$ such that $\be_{m-1}<r$ and $\ga_{1}>0$. Thus
$$s_D(b_1,\ldots,b_r,c_1,\ldots,c_s)=\sum_{\de\in \Delta_{r,s}} s_{\de},$$
where for $\de=(\be_1,\ldots,\be_m,\ga_1,\ldots,\ga_m)$ we denote by  $s_{\de}$ the following sum:
$$s_{\de}=s_{\de}(b_1,\ldots,b_r,c_1,\ldots,c_s)=\sum_{\si\in S_r,\;\tau\in S_s} b_{\si}^{\be}(1) c_{\tau}^{\ga}(1)\cdots  b_{\si}^{\be}(m) c_{\tau}^{\ga}(m).$$

If $\be_{i+1}-\be_i\geq d_1$ or $\ga_{i+1}-\ga_i\geq d_2$ for some $0\leq i<m$, then $s_{\de}=0$ in $V$, since $V_1$ and $V_2$ satisfy the symmetric identities of degrees $d_1$ and $d_2$, respectively.

Assume that $\be_{i+1}-\be_i< d_1$ and $\ga_{i+1}-\ga_i< d_2$ for all $0\leq i<m$.
Since
\begin{multline*}
D=(\be_1-\be_0)+(\be_2-\be_1)+\cdots+(\be_m-\be_{m-1})+{}\\
(\ga_1-\ga_0)+(\ga_2-\ga_1)+\cdots+(\ga_m-\ga_{m-1}),
\end{multline*}
we obtain $D\leq m(d_1+d_2-2)$. Hence  $m\geq kd_3+2$.

For a fixed $\de$ and some $i\geq 2$, $\si\in S_r$, $\tau\in S_s$, we write $w_{\si,\tau}^{\de}(i)$ for the product
$$b_{\si}^{\be}(i)c_{\tau}^{\ga}(i)\cdots b_{\si}^{\be}(i+k-1)  c_{\tau}^{\ga}(i+k-1).$$
Obviously, $b_{\si}^{\be}(1)$ is an element of $\hat{V}_1\cup\{1\}$, $b_{\si}^{\be}(2),\ldots,b_{\si}^{\be}(m)$ are elements of $\hat{V}_1$, $c_{\tau}^{\ga}(1),\ldots,c_{\tau}^{\ga}(m-1)$ are elements of $\hat{V}_2$ and $c_{\tau}^{\ga}(m)$ is an element of $\hat{V}_2\cup\{1\}$. Therefore, $w_{\si,\tau}^{\de}(i)\in (\hat{V}_1\hat{V}_2)^k$ in case  $1<i\leq m-k$.
The inequality $m\geq kd_3+2$ implies that we can rewrite $s_{\de}$ as follows:
\begin{multline*}
s_{\de}=\sum b_{\si}^{\be}(1) c_{\tau}^{\ga}(1)\times\\ 
\left(\sum w^{\de}_{\si,\tau}(2)\, w^{\de}_{\si,\tau}(k+2)\, w^{\de}_{\si,\tau}(2k+2) \cdots w^{\de}_{\si,\tau}((d_3-1)k+2)\right)\times\\
b_{\si}^{\be}(k d_3+2) c_{\tau}^{\ga}(k d_3+2)\cdots b_{\si}^{\be}(m) c_{\tau}^{\ga}(m),
\end{multline*}
where the first sum ranges over all $\si\in S_r$, $\tau\in S_s$ with $$\si(2)<\si(k+2) <\si(2k+2) <\cdots< \si((d_3-1)k+2)$$
and the second sum permutes the $w$'s. Since $(\hat{V}_1\hat{V}_2)^k$ satisfies the identity $s_{d_3}=0$, the second sum is zero and $s_{\de}$ is zero as well. Thus $V$ satisfies the identity $s_{D}(a_1,\ldots,a_D)=0$.
\end{proof}

Let us prove Theorem~\ref{theo_main}.

\begin{proof}
First, assume that $(R_1R_2)^k\subset R_1$. There exists a vector subspace $U_i$ of $R_i$ of dimension $t_i$ such that $R_i$ is the direct sum of the subspaces $L_i$ and $U_i$ ($i=1,2$). Let $l_{i1},l_{i2}\ldots$ be an $\FF$-basis for $L_i$ and $u_{i1},\ldots,u_{it_i}$ be an $\FF$-basis for $U_i$. Since the function $s_D$ is linear, without loss of generality we can assume that $a_1,\ldots,a_D$ are elements of the bases under consideration. Denote by $f_i$ the number of elements of the set $a_1,\ldots,a_D$ from the basis of $U_i$. If $f_i> t_i(p-1)$ for $i=1$ or $i=2$, then $s_D(a_1,\ldots,a_D)=0$ by Lemma~\ref{lemma1}.

Assume that $f_i\leq t_i(p-1)$ for $i=1,2$. Then $s_D(a_1,\ldots,a_D)$ is the sum of all elements of the following form:
$$y_0 z_1 y_1 \cdots z_q y_q,$$
where $y_0,\ldots,y_q$  are products of elements of the bases $l_{11},l_{12}\ldots$ and $l_{21},l_{22}\ldots$ of $L_1$ and $L_2$, respectively, and $z_1,\ldots,z_q$ are products of elements of the bases $u_{11},\ldots,u_{1t_1}$ and  $u_{21},\ldots,u_{2t_2}$ of $U_1$ and $U_2$, respectively. Here the only monomials that can be empty are $y_0$ and $y_q$. In the rest of the proof the degree of some product $a$ of elements of the set $\{y_0,\ldots,y_q,z_1,\ldots,z_q\}$ is the degree of $a$ as a monomial in $$l_{11},l_{12},\ldots,l_{21},l_{22},\ldots,u_{11},\ldots,u_{1t_1},u_{21},\ldots,u_{2t_2}. $$

For short, we write $\al=(t_1+t_2)(p-1)$ and $\beta=(d_1+d_2-2)(kd_1+2)$. Since the degree of $z_1\cdots z_q$ is $f_1+f_2\leq \al$, then $q\leq \al$. Therefore, there exists $j$ such that the degree of $y_{j}$ is greater than or equal to $\be$, since otherwise $D\leq (\be-1)(\al+1)+\al=\be(\al+1)-1$, a contradiction. Lemma~\ref{lemma_main} implies that $L_1+L_2$ satisfies the symmetric identity of degree $\beta$. It follows from the last two facts that $s_D(a_1,\ldots,s_D)=0$.

The proof in the case of $(R_1R_2)^k\subset R_2$ is the same.
\end{proof}

\section{Low bounds}

In this section we assume that the characteristic of $\FF$ is arbitrary. An identity $f$ is called \emph{multilinear} if $f=\sum_{\si\in S_n} \al_{\si} x_{\si(1)}\cdots x_{\si(n)}$, where $\al_{\si}\in \FF$. Consider an identity $f=\sum_i \al_i w_i$, where $w_i$ is a monomial in $x_1,x_2,\ldots$ and $\al_i\in \FF$. The identity $f$ is called \emph{homogeneous} if $\deg(w_i)=\deg(w_j)$ for all $i,j$. The degree of the identity $f$ is the maximal degree of the monomials $\{w_i\}$. The \emph{low degree} of the identity $f$ is the minimal degree of the  monomials $\{w_i\}$. Note that the low degree of an identity is always positive. It is well-known that if an algebra $R$ satisfies an identity of degree $d>0$, then $R$ satisfies a multilinear identity of degree $d$. Moreover, in the case of an infinite field if an algebra $A$ satisfies an identity of low degree $d$, then $A$ satisfies a multilinear identity of degree $d$.

\begin{example}
Let $\FF=\FF_2$ be the field of two elements and let an $\FF$-algebra $A$ be the quotient of the free associative algebra over the T-ideal generated by the identity $f(x)=x^2+x$. Then $A$ satisfies the multilinear identity $xy+yx$ of degree two, but does not satisfy any identity of degree one.
\end{example}

\begin{theo}\label{theo_low_bound} Let $\FF$ be a field of an arbitrary characteristic. Then for every identities $f_1$ and $f_2$ of low degrees $d_1>0$ and $d_2>0$, respectively, there exists an algebra $R=R_1+R_2$ such that
\begin{enumerate}
\item[$\bullet$] $R_i$ satisfies the identity $f_i$ for $i=1,2$ (and $R_i$ does not satisfy any identity of low degree less than $d_i$);

\item[$\bullet$] $R_2$ is a two-sided ideal in $R$;
\end{enumerate}
but $R$ does not satisfy any identity of degree less than $d_1 d_2$.
\end{theo}
\begin{proof}
Since the identity of nilpotency $x_1\cdots x_d=0$ implies any other identity of low degree $d$, it is enough to prove the theorem in the case of  $f_i=x_1\cdots x_{d_i}$  for $i=1,2$.

Let us construct the following example. Denote by $B$ the free associative algebra (without unity) freely generated by the elements $y_1,y_2,\dots,z_1,z_2,\dotsc$. Define the subalgebras $B_1$ and $B_2$ of $B$ as follows:
\begin{enumerate}
\item[$\bullet$] $B_1$ is generated by $y_1,y_2,\dotsc$;

\item[$\bullet$] $B_2$ is generated by products $u=u_{i_1}\cdots u_{i_n}$, where $u_i$ is $y_i$ or $z_i$ for $1\leq i\leq n$ and $u$ contains at least one element of the set $\{z_1,z_2,\ldots\}$.
\end{enumerate}
Denote by $R$ the quotient of $B$ over the two-sided ideal $I$ generated by all products $u_1\cdots u_{d_1}$ and $v_1\cdots v_{d_2}$, where $u_i\in B_1$ and $v_j\in B_2$ for all $i,j$. We write $R_1$ and $R_2$ for the images of $B_1$ and $B_2$ in $R$, respectively. Then $R=R_1\oplus R_2$ as vector spaces and $R_2$ is a two-sided ideal in $R$.

Assume that there exists an identity $f$ of degree $d=d_1d_2-1$ that holds in $R$. As we have already mentioned, without loss of generality we can assume that that $f$ is multilinear, i.e., $f$ is a sum of the monomials $\al_{\si} x_{\si(1)} \cdots x_{\si(d)}$, where $\si\in S_d$ and $\al_{\si}\in\FF$. Moreover, we can assume that  $\al_{\rm id}=1$ for the trivial permutation ${\rm id}\in S_d$. In other words, $f=x_1\cdots x_d+h$, where the monomial $x_1\cdots x_d$ does not appear in the polynomial $h$. Define the elements $a_1,\ldots,a_d$ of $R$ as follows:
$$a_i=\left\{
\begin{array}{rl}
z_j & \text{ if } i=j d_1 \text{ for some }j>0,\\
y_i & \text{ otherwise. } \\
\end{array}
\right.$$
Then 
\begin{multline*}
f(a_1,\ldots,a_d)=h(a_1,\ldots,a_d)+{}\\
+\underbrace{y_1\cdots y_{d_1-1}\, z_1}_{d_1} \cdots \underbrace{y_{d_1(d_2-2)+1} \cdots y_{d_1(d_2-1)-1} \,z_{d_2-1}}_{d_1}\, y_{d_1(d_2-1)+1}\cdots y_{d_1d_2-1}. 
\end{multline*}
Note that the last monomial of $f(a_1,\ldots,a_d)$ does not appear as a monomial in any element of the ideal $I$. Then it is not difficult to see that $f(a_1,\ldots,a_d)\neq0$ in $R$; a contradiction. The theorem has been proved.
\end{proof}

Note that the semigroup algebra $\FF[S]$ for the semigroup $S$ from the last example from A.~Salwa's paper~\cite{Salwa96} has the same properties as that presented in the proof of Theorem~\ref{theo_low_bound}. The proof of Theorem~\ref{theo_low_bound} implies the following remark.

\begin{remark}
The statement of Theorem~\ref{theo_low_bound} also holds for algebras with unity. Namely,
for every $f_1,f_2$ from $\FF\LA X\RA$ with $f_i(1,\ldots,1)=0$ for $i=1,2$ there exists an algebra $R$ with unity $1$ and with subalgebras $R_1,R_2$ such that $R_i$ contains $1$ ($i=1,2$) and the condition of Theorem~\ref{theo_low_bound} holds for $R=R_1+R_2$.
\end{remark}

\section*{Acknowledgements}
The second author was supported by FAPESP No.~2013/15539-2 and RFFI 14-01-00068. He is grateful for this support.



\begin{thebibliography}{99}
\bibitem{BahturinGiambruno94} Yu. Bahturin, A. Giambruno, \emph{Identities of sums of commutative subalgebras}, Rend. Cir. Mat. Palermo (2), {\bf 43} (1994), 250--258.

\bibitem{BeidarMikhalev95} K.I. Beidar, A.V. Mikhalev, \emph{Generalised polynomial identities and rings which are sums of two subrings}, Algebra i Logica, {\bf 34} (1995), 3--11 (Russian). English translation: Algebra and Logic, {\bf 34} (1995), 1--5.

\bibitem{Bokut76} L.A. Bokut', \emph{Embeddings in simple associative algebras}, Algebra i Logica, {\bf 15} (1976), 117--142 (Russian). English translation: Algebra and Logic, {\bf 15} (1976), 73--90.

\bibitem{FelzenszwalbGiambrunoLeal03} B. Felzenszwalb, A. Giambruno,  G. Leal, \emph{On rings which are sums of two PI-subrings: a combinatorial approach}, Pacific J. Math., {\bf 209} (2003), no. 1, 17--31.

\bibitem{Goto62} M. Goto, \emph{Note on a characterization of solvable Lie algebras}, J. Sci. Hiroshima Univ. Ser. A-I {\bf 26} (1962), 1--2.

\bibitem{Kegel62} O.H. Kegel, \emph{Zur Nilpotenz gewisser assoziativer Ringe}, Math. Ann., {\bf 149} (1962/63), 258--260.

\bibitem{Kemer93} A.R. Kemer, \emph{The standard identity in charcteristic $p$: a conjecture of I.B. Volichenko}, Israel J. Math., {\bf 81} (1993), 343--355.

\bibitem{KepczykPuczylowski96} M. K\c{e}pczyk, E.R. Puczy\l{}owski, \emph{On radicals of rings which are sums of two subrings}, Arch. Math., {\bf 66} (1996), 8--12.

\bibitem{KepczykPuczylowski98} M. K\c{e}pczyk, E.R. Puczy\l{}owski, \emph{Rings which are sums of two subrings}, J. Pure Applied Algebra, {\bf 133} (1998), 151--162.

\bibitem{KepczykPuczylowski01} M. K\c{e}pczyk, E.R. Puczy\l{}owski, \emph{Rings which are sums of two subrings satisfying a polynomial identity}, Commun. Algebra, {\bf 29} (2001), 2059--2065.

\bibitem{Kepczyk08} M. K\c{e}pczyk, \emph{On algebras that are sums of two subalgebras satisfying certain polynomial identities}, Publ. Math. Debrecen, {\bf 72/3} (2008), 257--267.

\bibitem{Kepczyk15} M. K\c{e}pczyk, \emph{Note on algebras which are sums of two PI subalgebras}, J. Algebra Appl., {\bf 14}, (2015) no.~10, 1550149.

\bibitem{Kolman65} B. Kolman, \emph{Semi-modular Lie algebras}, J. Sci. Hiroshima Univ. Ser. A-I, {\bf 29} (1965), 149--163.

\bibitem{Panyukov90} V.V. Panyukov, \emph{On the solvability of Lie algebras of positive characteristic}, Russ. Math. Surv., {\bf 45} (1990), No.4, 181--182.

\bibitem{Petravchuk88} A.P. Petravchuk, \emph{Lie algebras which can be decomposed into the sum of an abelian subalgebra and a nilpotent subalgebra}, Ukrain. Math. J., {\bf 40} (1988), 331--334.

\bibitem{Pchelintsev85} S.V. Pchelintsev, \emph{Kegel theorem for alternative algebras}, Sibirsk. Mat. Zh., {\bf 26} (1985), 195--196.

\bibitem{Regev71} A. Regev, \emph{Existence of polynomial identities in $A\otimes_F B$}, Bull. Amer. Math. Soc., {\bf 77} (1971), 1067--1069.

\bibitem{Regev72} A. Regev, \emph{Existence of identities in $A\otimes B$}, Israel J. Math., {\bf 11} (1972), 131--152.

\bibitem{Rowen76} L.H. Rowen, \emph{Generalized polynomial identities II}, J. Algebra, {\bf 38} (1976), 380--392.

\bibitem{Salwa96} A. Salwa, \emph{Rings that are sums of two locally nilpotent subrings}, Commun. Algebra, {\bf 24} (1996), no. 12, 3921--3931.

\bibitem{Zusmanovich91} P. Zusmanovich, \emph{A Lie algebra that can be written as a sum of two nilpotent subalgebras is solvable}, Mathematical Notes, {\bf 50}  (1991), 909--912.
\end{thebibliography}
\end{document}